\documentclass[10pt, a4paper]{article}
\usepackage{amscd}
\usepackage{amsfonts}
\usepackage{caption}
\usepackage{color}
\begin{document}

\newtheorem{theorem}{Theorem}
\newtheorem{corollary}[theorem]{Corollary}
\newtheorem{definition}[theorem]{Definition}
\newtheorem{lemma}[theorem]{Lemma}
\newtheorem{proposition}[theorem]{Proposition}
\newtheorem{remark}[theorem]{Remark}
\newtheorem{example}[theorem]{Example}
\newtheorem{notation}[theorem]{Notation}
% fim da demonstracao
\def\Qed{\hfill\raisebox{.6ex}{\framebox[2.5mm]{}}\\[.15in]}

\title{Some bidouble planes with $p_g = q = 0$ and $4\leq K^2\leq 7$}

\author{Carlos Rito}

\date{}
\pagestyle{myheadings}
\maketitle
\setcounter{page}{1}

\begin{abstract}
We give a list of possibilities for surfaces of general type with $p_g=0$ having an involution $i$ such that the bicanonical map of $S$ is not composed with $i$ and $S/i$ is not rational.
Some examples with $K^2=4,\ldots,7$ are constructed as double coverings of an Enriques surface. These surfaces have a description as bidouble coverings of the plane.

\end{abstract}

\section{Introduction}

Smooth minimal surfaces $S$ of general type with $p_g=0$ have been studied by several authors in the last years, but a classification is still missing.
We refer the survey \cite{BCP} for information on these surfaces.

The case where $S$ has an involution $i$ such that the bicanonical map of $S$ is composed with $i$ has been considered by Mendes Lopes and Pardini upon their study on surfaces with $p_g=0$ having non-birational bicanonical map (\cite{MP1}, \cite{MP2}, \cite{MP3}, \cite{MP4}). Due to the work of Du Val (\cite{Du}), Xiao Gang (\cite{Xi}) and Borrelli (\cite{Bo}), there is a classification for the case where the bicanonical map is of degree 2 onto a rational surface.

In this paper we consider the case where the bicanonical map is not composed with $i$ and $S/i$ is not a rational surface. Since $p_g(S/i)=q(S/i)=0,$ it is natural to expect to have a short list of possibilities. We confirm this with a structure theorem and we give some examples. We do not claim that these surfaces are new; our aim is to focus on the involutions. It was easier for us to construct these examples than to identify involutions on the already long list of known examples of surfaces with $p_g=0$, in particular product-quotient surfaces.  

These examples provide the existence of the cases with $K^2=5,6,7$ and $S/i$ birational to an Enriques surface (cf. arXiv:1003.3595v2, \cite{LS}). More recently Yifan Chen (\cite{Ch}) has shown that our example with $K^2=7$ is an Inoue surface (\cite{In}).

The paper is organized as follows.
First we note that some modifications to \cite[Theorems $7,$ $8$ and $9$]{Ri1} give a list of possibilities for the branch curve in the quotient surface $S/i$.
Then we construct some examples with $K_S^2=4,\ldots,7$ as double covers of an Enriques surface, which in turn is obtained as a quotient of a Kummer surface.
Finally we describe these surfaces as bidouble covers of the plane and we give some other bidouble plane examples.
In all cases $S$ has another involution $j$ such that $S/j$ is rational and the bicanonical map of $S$ is composed with $j$.

\bigskip
\noindent{\bf Notation}

We work over the complex numbers; all varieties are assumed to be projective algebraic.
An {\em involution} of a surface $S$ is an
automorphism of $S$ of order 2. We say that a map is {\em composed with
an involution} $i$ of $S$ if it factors through the double cover $S\rightarrow
S/i.$
A {\em $(-2)$-curve} or {\em nodal curve} $N$ on a surface is a curve isomorphic to $\mathbb P^1$ such that $N^2=-2$.
An $(m_1,m_2,\ldots)$-point of a curve, or point of type $(m_1,m_2,\ldots),$ is a singular point of multiplicity $m_1,$ which resolves to a point of multiplicity $m_2$ after one blow-up, etc.
The rest of the notation is standard in Algebraic Geometry.\\

\bigskip
\noindent{\bf Acknowledgements}

The author wishes to thank Margarida Mendes Lopes for all the support.
He is a member of the Mathematics Center of the Universidade de Tr\'as-os-Montes e Alto Douro and is a collaborator of the Center for Mathematical Analysis, Geometry and Dynamical Systems of Instituto Superior T\'ecnico, Universidade T\' ecnica de Lisboa.

This research was partially supported by the Funda\c c\~ao para a Ci\^encia e a Tecnologia (Portugal) through Projects PEst-OE/MAT/UI4080/2011 and\newline PTDC/MAT/099275/2008.

\section{General facts on involutions}\label{GenFacts}

The following is according to \cite{CM}.\\
Let $S$ be a smooth minimal surface of general type with an
involution $i.$ Since $S$ is minimal of general type, this
involution is biregular. The fixed locus of $i$ is the union of a
smooth curve $R''$ (possibly empty) and of $t\geq 0$ isolated points
$P_1,\ldots,P_t.$ Let $S/i$ be the quotient of $S$ by $i$ and
$p:S\rightarrow S/i$ be the projection onto the quotient. The
surface $S/i$ has nodes at the points $Q_i:=p(P_i),$ $i=1,\ldots,t,$
and is smooth elsewhere. If $R''\not=\emptyset,$ the image via $p$
of $R''$ is a smooth curve $B''$ not containing the singular points
$Q_i,$ $i=1,\ldots,t.$ Let now $h:V\rightarrow S$ be the blow-up of
$S$ at $P_1,\ldots,P_t$ and set $R'=h^*(R'').$ The involution $i$
induces a biregular involution $\widetilde{i}$ on $V$ whose fixed
locus is $R:=R'+\sum_1^t h^{-1}(P_i).$ The quotient
$W:=V/\widetilde{i}$ is smooth and one has a commutative diagram:
$$
\begin{CD}\ V@>h>>S\\ @V\pi VV  @VV p V\\ W@>g >> S/i
\end{CD}
$$
where $\pi:V\rightarrow W$ is the projection onto the quotient and
$g:W\rightarrow S/i$ is the minimal desingularization map. Notice
that $$A_i:=g^{-1}(Q_i),\ \ i=1,\ldots,t,$$ are $(-2)$-curves and
$\pi^*(A_i)=2\cdot h^{-1}(P_i).$

Set $B':=g^*(B'').$ Since $\pi$ is
a double cover with branch locus $B'+\sum_1^t A_i,$ it is determined
by a line bundle $L$ on $W$ such that $$2L\equiv B:=B'+\sum_1^t A_i.$$

\begin{proposition}[\cite{CM}, \cite{CCM}]
The bicanonical map of $S$ (given by $|2K_S|$) is composed with $i$ if and only if $h^0(W,\mathcal{O}_W(2K_W+L))=0$.
\end{proposition}

\section{List of possibilities}

Let $P$ be a minimal model of the resolution $W$ of $S/i,$ let $\rho:W\rightarrow P$ be the corresponding projection and denote by $\overline{B}$ the projection $\rho(B).$

\begin{theorem}\label{thm}
{\em ({\rm cf.} \cite{Ri1})}
Let $S$ be a smooth minimal surface of general type with $p_g=0$ having an involution $i$ such that the bicanonical map of $S$ is not composed with $i$ and $S/i$ is not rational.

With the previous notation, one of the following holds:
\begin{description}
  \item[a)] $P$ is an Enriques surface and:

    $\cdot$ $\overline B^2=0,$ $t-2=K_S^2\in\{2,\ldots,7\},$ $\overline B$ has a $(3,3)$-point or a $4$-uple point and at most one double point.

  \item[b)] ${\rm Kod}(P)=1$ and:

    $\cdot$ $K_P\overline B=2,$ $\overline B^2=-12,$ $t-2=K_S^2\in\{2,\ldots,8\},$ $\overline B$ has at most two double points, or

    $\cdot$ $K_P\overline B=4,$ $\overline B^2=-16,$ $t=K_S^2\in\{4,\ldots,8\},$ $\overline B$ is smooth.

  \item[c)] ${\rm Kod}(P)=2$ and:

    $\cdot$ $K_S^2=2K_P^2,$ $K_P^2=1,\ldots,4,$ $\overline B$ is a disjoint union of four $(-2)$-curves, or

    $\cdot$ $K_P\overline B=2,$ $K_P^2=1,$ $\overline B^2=-12,$ $t=K_S^2\in\{4,\ldots,8\},$ $\overline B$ has at most one double point, or

    $\cdot$ $K_P\overline B=2,$ $K_P^2=2,$ $\overline B^2=-12,$ $t+2=K_S^2\in\{6,7,8\},$ $\overline B$ is smooth.
\end{description}
Moreover there are examples for {\rm a), b)} and {\rm c)}.
\end{theorem}
{\bf Proof :} This follows from the proof of \cite[Theorems $7,$ $8$ and $9$]{Ri1} taking in account that:
\begin{enumerate}
\item[$\cdot$] $p_g(P)=q(P)=0$ (because $p_g(P)\leq p_g(S),$ $q(P)\leq q(S)$);
\item[$\cdot$] $h^0(W,\mathcal{O}_W(2K_W+L))\leq\frac{1}{2}K_W^2+2$ (see \cite[Proposition $4$, b)]{Ri1});
\item[$\cdot$] $K_S^2\ne 9$ (see \cite[Theorem $4.3$]{DMP});
\item[$\cdot$] We can have $\overline {B'}^2>0$ (unlike the case $p_g=q=1$).

Examples for a) and b) are given below. Rebecca Barlow (\cite{Ba}) has constructed a surface of general type with $p_g=0$ and $K^2=1$ containing an even set of four disjoint $(-2)$-curves. This gives an example for c).

\end{enumerate}

\section{Examples}

\subsection{$S/i$ birational to an Enriques surface}\label{Enr}

Consider the involution of $\mathbb P^1\times \mathbb P^1$ $$j:[x:y\ ,\ a:b]\mapsto [y:x\ ,\ b:a]$$ and denote by $f,g$ the projections onto the first and second factors, respectively.
Let $F_1,\ldots,F_4,$ $G_1,\ldots,G_4$ be fibres of $f,$ $g$ such that $$C:=F_1+\cdots F_4+G_1+\cdots+G_4$$ is preserved by $j$ and does not contain the fixed points $[1:\pm 1,1:\pm 1]$ of $j.$

Let $$\pi:Q\rightarrow \mathbb P^1\times \mathbb P^1$$ be the double cover with branch locus $C$ and let $k$ be the corresponding involution.
It is well known that $Q$ is a Kummer surface and $$E:=Q/{k\circ j}$$ is an Enriques surface with $8$ nodes.

\subsubsection{$\overline B$ with a $4$-uple point}\label{Enr4}

Let $D\subset\mathbb P^1\times \mathbb P^1$ be a generic curve of bi-degree $(1,2)$ tangent to $C$ at smooth points $p_1, p_2$ of $C$ such that $p_2=j(p_1).$ The pullback $\pi^*(D+j(D))\subset Q$ is a reduced curve with two $4$-uple points, corresponding to the $(2,2)$-points of $D+j(D)$ (which are tangent to the branch curve $C$).
These points are identified by the involution $k\circ j,$ thus the projection of $\pi^*(D+j(D))$ into $E$ is a reduced curve $\overline{\overline {B'}}$ with one $4$-uple point.

Now let $\widetilde{E}$ be the minimal smooth resolution of the Enriques surface $E,$ $A_1,\ldots,A_8\subset\widetilde{E}$ be the nodal curves corresponding to the nodes of $E$ and $\overline {B'}\subset\widetilde E$ be the strict transform of $\overline{\overline {B'}}.$
If $D$ does not contain one of the $16$ double points of $C,$
the divisor $$\overline B:=\overline{B'}+\sum_1^8A_i$$ is reduced, divisible by $2$ in the Picard group and satisfies $\overline B^2=0.$
Let $S$ be the smooth minimal model of the double cover of $\widetilde E$ ramified over $\overline B$.
One can show that $S$ is a surface of general type with $p_g=0$ and $K_S^2=6.$
Moreover, $D$ can be chosen through one or two double points of $C.$ This provides examples with $K_S^2=5$ or $4$,
corresponding to branch curves $$\overline{B'}+\sum_1^7A_i\ \ \ \ {\rm or}\ \ \ \ \overline{B'}+\sum_1^6A_i.$$

\subsubsection{$\overline B$ with a $(3,3)$-point}\label{765}

Let $D_1$ be a curve of bi-degree $(0,1)$ through $p,$ $D_2$ be a general curve of bi-degree $(1,1)$ through $p$ and $j(p)$ and set $D:=D_1+D_2.$ Then $D+j(D)$ is a reduced curve with triple points at $p$ and $j(p).$ Now we proceed as in Section \ref{Enr4}. In this case the branch curve $\overline B\subset\widetilde E$ has a $(3,3)$-point instead of a $4$-uple point. This gives an example of a surface of general type $S$ with $p_g=0$ and $K^2=7$ (notice that the resolution of the $(3,3)$-point gives rise to an additional nodal curve in the branch locus). As above, $D$ can be chosen containing one or two double points of $C,$ providing examples with $K_S^2=6$ or $5$.

\subsection{Bidouble plane description}\label{Enr33}

Here we obtain the examples of Section \ref{Enr} as bidouble covers of the plane.

\subsubsection{Construction}\label{construction}

Let $T_1,\ldots,T_4\subset\mathbb P^2$ be distinct lines through a point $p$ and $C_1,$ $C_2$ be distinct smooth conics tangent to $T_1,$ $T_2$ at points $p_1,p_2\ne p,$ respectively.
The smooth minimal model $\widetilde E$ of the double cover of $\mathbb P^2$ with branch locus $T_1+\ldots+T_4+C_1+C_2$ is an Enriques surface with $8$ disjoint nodal curves $A_1,\ldots,A_8,$ which correspond to the $8$ double points of $$G:=T_3+T_4+C_1+C_2.$$

Now let $p_3$ be a generic point in $T_3$ and consider the pencil $l$ generated by $2H_i+T_i,$ $i=1,2,3,$ where $H_i$ is a conic through $p_i$ tangent to $T_j,T_k$ at $p_j,p_k,$ for each permutation $(i,j,k)$ of $(1,2,3).$
Let $L$ be a generic element of $l.$
Notice that the quintic curve $L$ contains $p,$ it has a $(2,2)$-point at $p_i$ and the intersection number of $L$ and $T_i$ at $p_i$ is $4,$ $i=1,2,3.$

The strict transform of $L$ in $\widetilde E$ is a reduced curve $\overline{B'}$ with a $4$-uple point (at the pullback of $p_3$) such that the divisor $$\overline B:=\overline{B'}+\sum_1^8A_i$$
is reduced, satisfies $\overline B^2=0$ and is divisible by $2$ in the Picard group (because $L+T_1$ is divisible by $2$).
Let $S$ be the smooth minimal model of the double cover of $\widetilde E$ ramified over $\overline B$.
One can verify that $K_S^2=6.$
As in Section \ref{Enr4}, choosing $L$ through $1$ or $2$ double points of $T_3+T_4+C_1+C_2$ one obtains examples with $K_S^2=5$ or $4,$ respectively.

To obtain a branch curve $\overline B\subset\widetilde E$ with a $(3,3)$-point as in Section \ref{Enr33}, it suffices to change the $(2,2)$-point of the quintic $L$ at $p_3$ to an ordinary triple point. In this case $L$ is the union of a conic through $p_3$ with a cubic having a double point at $p_3.$ Choosing $C_1$ and $C_2$ so that $L$ passes through $0,$ $1$ or $2$ double points of $T_3+T_4+C_1+C_2$ one obtains examples with $K_S^2=7,$ $6$ or $5.$

\subsubsection{Involutions on $S$}

We refer \cite{Ca} or \cite{Pa} for information on bidouble covers.

Each surface $S$ constructed in Section \ref{construction} is the
smooth minimal model of the bidouble cover of $\mathbb P^2$
determined by the divisors
$$
\begin{array}{l}
D_1:=L,\\
D_2:=T_1+C_1+C_2,\\
D_3:=T_2+T_3+T_4.
\end{array}
$$

Let $i_g$ be the involution of $S$ corresponding to $D_j+D_k,$ for
each permutation $(g,j,k)$ of $(1,2,3).$ We have
that $S/{i_1}$ is birational to an Enriques surface, $S/{i_3}$ is a rational surface and
the bicanonical map of $S$ is not composed with $i_1, i_2$ and is composed with $i_3.$
The surface $S/{i_2}$ is birational to an Enriques surface in the cases with $K_S^2=6,5,4$ given in Section \ref{Enr4} and is a rational surface in the cases with $K_S^2=7,6,5$ given in Section \ref{765}.
Moreover, $S$ has an hyperelliptic fibration of genus $3$.

We omit the proof for these facts: it is similar to the one given in \cite{Ri3} for an example with $K_S^2=3.$

\subsection{More bidouble planes}

In the examples above, $S/{i_1}$ is birational to an Enriques
surface with $8$ disjoint $(-2)$-curves, corresponding to the $8$
nodes of the sextic $G=T_3+T_4+C_1+C_2,$ which contains $2$ lines.
%is a double plane with branch curve the union of two lines $T_1,T_2$ with a sextic $G$ which contains two other lines.
Now we give examples with $G$ containing only one line and with $G$ without lines.

\subsubsection{$G$ with one line, $K_S^2=4,5,6$}

Let $T_1,T_2,T_3$ and $L$ be as in Section \ref{construction} and $p_4$ be a smooth point of $L.$ There exists a plane curve $J$ of degree $5$ through $p$ with $(2,2)$-points tangent to $T_1,T_2,L$ at $p_1,p_2,p_4,$ respectively (notice that we are imposing $19$ conditions to a linear system of dimension $20$; such a curve can be easily computed using the Magma function LinSys given in \cite{Ri2}).

Let $S$ be the smooth minimal model of the bidouble cover of
$\mathbb P^2$ determined by the divisors
$$
\begin{array}{l}
D_1:=L,\\
D_2:=T_3,\\
D_3:=T_1+T_2+J.
\end{array}
$$
 Notice that the double plane with branch locus $D_2+D_3$ is an
Enriques surface $E$ with $6$ disjoint nodal curves
$A_1,\ldots,A_6$ (two of them are contained in the pullback of
$p_4$) and that the strict transform $\widehat L$ of $L$ in $E$
has a $4$-uple point at the pullback of $p_3.$ Moreover, the
divisor $\overline B:=\widehat L+\sum_1^6 A_i$ satisfies
$\overline B^2=0$ and is even (because $L+T_3$ is even).
This gives an example for Theorem \ref{thm}, a) with $K_S^2=4.$

To obtain an example with $K_S^2=5$ it suffices to choose the quintic $J$ with a triple point at $p_4$ instead of a $(2,2)$-point. In this case $J$ is the union of a conic with a singular cubic. Here the Enriques surface contains $7$ disjoint nodal curves, three of them contained in the pullback of $p_4.$

Finally, choosing $L$ with a triple point at $p_3$ one obtains $\widehat L\subset E$ with a $(3,3)$-point.
This gives examples for Theorem \ref{thm}, a) with $K_S^2=5,6.$

\subsubsection{$G$ without lines, $K_S^2=4$}

Consider, in affine plane, the points $p_0,\ldots,p_5$ with coordinates $(0,0),$ $(1,1),$ $(-1,1),$ $(2,3),$ $(-2,3),$ $(0,5),$ respectively, and let $T_{ij}$ be the line through $p_i,p_j.$
Let $C_1$ be the conic tangent to $T_{01},T_{02}$ at $p_1,p_2$ which contains $p_5$ and let $C_2$ be the conic tangent to $T_{01},T_{02}$ at $p_1,p_2$ which contains $p_3,p_4.$
Let $l$ be the linear system generated by $T_{01}+T_{02}+2T_{34}$ and $T_{03}+T_{04}+C_2.$

The element $Q$ of $l$ through $p_5$ is an irreducible quartic curve with double points at $p_0,p_3,p_4$ and tangent to $T_{01},T_{02}$ at $p_1,p_2$. Moreover, because of the symmetry with respect to $T_{05},$ the line $H$ tangent to $Q$ at $p_5$ is horizontal.

There is a cubic $F$ through $p_0,p_3,p_4$ tangent to $T_{01},$ $T_{02},$ $H$ at $p_1,$ $p_2,$ $p_5,$ respectively (notice that we are imposing $9$ conditions to a linear system of dimension $9$).
One can verify that $F$ contains no line, thus it is irreducible.

The surface $S$ is the smooth minimal model of the bidouble cover of $\mathbb P^2$ determined by the divisors
$$
\begin{array}{l}
D_1:=C_1+F,\\
D_2:=T_{01},\\
D_3:=T_{02}+C_2+Q.
\end{array}
$$
Notice that the double plane with branch locus $D_2+D_3$ is an
Enriques surface $E$ with $6$ disjoint nodal curves (contained in
the pullback of the triple points $p_3,p_4$ of $D_3$) and that the
strict transform of $D_1$ in $E$ has a $4$-uple point (at the
pullback of $p_5$).
The double plane with branch locus $D_1+D_3$ is a surface with Kodaira dimension 1. This gives an example for Theorem \ref{thm}, a), b) with $K_S^2=4.$

\bibliography{ReferencesRito}

\bigskip
\bigskip

\noindent Carlos Rito
\\ Departamento de Matem\' atica
\\ Universidade de Tr\' as-os-Montes e Alto Douro
\\ 5001-801 Vila Real
\\ Portugal
\\\\
\noindent {\it e-mail:} crito@utad.pt

\end{document}